\documentclass[reqno]{amsart}
\usepackage{hyperref}
\usepackage[comma,sort&compress]{natbib}
\usepackage{fancyhdr,latexsym,amsmath,amsfonts,amssymb,amsbsy,amsthm,url,bbm}
\allowdisplaybreaks
\newtheorem{theorem}{Theorem}[section]
\newtheorem{lemma}{Lemma}[section]
\newtheorem{proposition}{Proposition}[section]

\usepackage{rotating}
\allowdisplaybreaks
\usepackage{enumerate}
\theoremstyle{definition}
\newtheorem{definition}{Definition}[section]

\theoremstyle{remark}
\newtheorem{remark}{Remark}[section]

\numberwithin{equation}{section}
\DeclareMathOperator{\E}{E}

\DeclareMathOperator{\prob}{P}

\DeclareMathOperator{\e}{e}
\newcommand{\comment}[1]{}
\newcommand{\Egamma}{\mathsf{E}^{(\gamma)}}
\newcommand{\Ebar}{\overline{\mathsf{E}}^{(\gamma)}}
\newcommand{\vaguec}{\stackrel{\mathrm{v}}{\rightarrow}}

\bibfont{\footnotesize}

\begin{document}
\bibliographystyle{abbrvnat}
\title{Products in Conditional Extreme Value Model}

\author[R. S. Hazra]{Rajat Subhra Hazra}
\address{Rajat Subhra Hazra\\Statistics and Mathematics Unit\\
Indian Statistical Institute\\ 203 B.T.~Road\\ Kolkata 700108\\
India} \email{rajat\_r@isical.ac.in}

\author[K. Maulik]{Krishanu Maulik}
\address{Krishanu Maulik\\ Statistics and Mathematics Unit\\ Indian
Statistical Institute\\ 203 B.T.~Road\\ Kolkata 700108\\ India}
\email{krishanu@isical.ac.in}


\begin{abstract}
The classical multivariate extreme value theory tries to capture the extremal dependence between the components under a multivariate domain of attraction condition and it requires each of the components to be in the domain of attraction of a univariate extreme value distribution as well. The multivariate extreme value (MEV) model has a rich theory but has some limitations as it fails to capture the dependence structure in presence of asymptotic independence. A different approach to MEV was given by \cite{Heffernan:Tawn:2004}, where they examined MEV distributions by conditioning on one of the components to be extreme. Here we assume one of the components to be in Frech\'et or Weibull domain of attraction and study the behavior of the product of the components under this conditional extreme value model.
\end{abstract}
\maketitle

\begin{section}{Introduction} \label{chap3: sec: intro}
The classical multivariate extreme value theory tries to capture the extremal dependence between the components under a multivariate domain of attraction condition and it requires each of the components to be in domain of attraction of a univariate extreme value distribution. The multivariate extreme value theory has a rich theory but has some limitations as it fails to capture the dependence structure. The concept of tail dependence is an alternative way of detecting this dependence. The concept was first proposed by \cite{ledford:tawn:1996,ledford:tawn:1997} and then elaborated upon by \cite{Resnick:2002, Maulik:Resnick:2004}. A different approach towards modeling multivariate extreme value distributions was given by \cite{Heffernan:Tawn:2004} by conditioning on one of the components to be  extreme. Further properties of this conditional model were subsequently studied by \cite{heffernan:resnick:2007,Das:Resnick:2008}.

An important dependence structure in multivariate extreme value theory is that of asymptotic independence. The joint distribution of two random variable is asymptotically independent if the the nondegenerate limit of suitably centered and scaled coordinate wise partial maximums is a product measure. One of the limitations of the asymptotic independence model is that it is too large a class to conclude anything interesting, for example, about product of two random variables. In another approach, a smaller class was considered by \cite{Maulik:Resnick:Rootzen:2002}. They assumed that $(X,Y)$ satisfy the following vague convergence:
\begin{equation}\label{asymptotic-ind-strong}
t\prob\left[ \left(\frac{X}{a(t)},Y\right)\in\cdot\right]\vaguec(\nu\times H)(\cdot) \text{  on  }M_+((0,\infty]\times[0,\infty]),
\end{equation}
where $M_+((0,\infty]\times[0,\infty])$ denotes the space of nonnegative Radon measures on $(0,\infty] \times [0,\infty]$ and $\nu(x,\infty]=x^{-\alpha}$, for some $\alpha>0$ and $H$ is a probability distribution supported on $(0,\infty)$. The tail behavior of the product $XY$  under the assumption~\eqref{asymptotic-ind-strong} and some further moment conditions was obtained by \cite{Maulik:Resnick:Rootzen:2002}. The conditional model can be viewed as an extension of the above model. Under the conditional model, the limit of the vague convergence need  not be a product measure and it happens in the space $M_+\left([-\infty,\infty]\times \Ebar\right)$ with $\gamma\in \mathbb{R}$, where $\Ebar$ is the right closure of the set $\{x\in\mathbb{R}:1+\gamma x>0\}$. (See Section~\ref{sec: notation} for details.) In this article we mainly focus on the product behavior when the limiting measure in the conditional model is not of the product form.

Products of random variables and their domain of attraction are important theoretical issues which have a lot of applications ranging from Internet traffic to insurance models. We study the product of two random variables whose joint distribution satisfies the conditional extreme value model. In particular we try to see the  role of regular variation in determining the behavior of the product. When $\gamma>0$, then it is easy to describe the behavior of the product under certain conditions. However, when $\gamma<0$, some complications arise due to the presence of finite upper end point. We remark that in this article we do not deal with the other important case of $\gamma =0$. Like in the case of Gumbel domain of attraction for maximum of i.i.d.\ observations, the case $\gamma=0$ will require more careful and detailed analysis.

In Section~\ref{sec: notation} we briefly describe the conditional extreme value model and state some properties and nomenclature which we use throughout this article. In Section~\ref{sec: overview} we provide an overview of our results presented in the later sections. In Section~\ref{Transformations} we reduce the conditional extreme value model defined in Section~\ref{sec: notation} to simpler forms in special cases. In Section~\ref{mainresults} we describe the behavior of the product of random variables following conditional model under some appropriate conditions. In Section~\ref{section:remarks on moments} we make some remarks on the assumptions used in the results. In the final Section~\ref{section:example} we present an example of a conditional model and also look at the product behavior.

\end{section}

\begin{section}{Conditional extreme value model} \label{sec: notation}
In this section, we provide the notations used in this article and the basic model.

If $\mathsf{S}$ is a topological space with $\mathcal{S}$ being its $\sigma$-field, then for non-negative Radon measures $\mu_t$, for $t>0$, and $\mu$ on $(\mathsf S, \mathcal S)$, we say $\mu_t$ converges vaguely to $\mu$ and denote it by $\mu_t \vaguec \mu$, if for all relatively compact sets $C$, which are also $\mu$-continuity sets of $\mu$, that is, $\mu( \partial C ) = 0$, we have $\mu_t (C) \to \mu(C)$, as $t \to \infty$. By $M_+(\mathsf S)$ we shall mean the space of all Radon measures on $\mathsf{S}$ endowed with the topology of vague convergence.

\begin{definition} \label{def: reg var}
A measurable function $f:\mathbb{R}_{+}\rightarrow \mathbb{R}_{+}$ is called regularly varying at infinity with index $\alpha$, (we write $f \in RV_{\alpha}$), if for all $t>0$, $$\lim_{x\rightarrow\infty}\frac{f(tx)}{f(x)} = t^\alpha.$$
If $\alpha = 0$, $f$ is called slowly varying. 
\end{definition}
\begin{definition} \label{def: reg var rv}
We say that a random variable $X$, with distribution function $F$, has a regularly varying tail of index $-\alpha$,  $\alpha \ge 0$, if the tail of its distribution function $\overline F(\cdot) := \prob[X>\cdot]\in RV_{-\alpha}$. 
\end{definition}
Definition~\ref{def: reg var rv} is equivalent to the existence of a positive function $a\in RV_{1/\alpha}$ such that $t\prob[X/a(t)\in \cdot]$ has a vague limit in $M_+((0,\infty])$, where the limit is a nondegenerate Radon measure. The limiting measure necessarily takes values $cx^{-\alpha}$ on set $(x,\infty]$.

Let $\Egamma$ be the interval $\{x \in \mathbb{R}: 1+\gamma x > 0\}$ and $\Ebar$ be its right closure in the extended real line $[-\infty, \infty]$. Thus, we have
$$\Egamma =
\begin{cases}
  (-1/\gamma, \infty), &\text{if $\gamma>0$},\\
  (-\infty, \infty), &\text{if $\gamma=0$},\\
  (-\infty, -1/\gamma), &\text{if $\gamma<0$},\\
\end{cases}
\qquad \text{and} \qquad
\Ebar =
\begin{cases}
  (-1/\gamma, \infty], &\text{if $\gamma>0$},\\
  (-\infty, \infty], &\text{if $\gamma=0$},\\
  (-\infty, -1/\gamma], &\text{if $\gamma<0$}.\\
\end{cases}$$

For any $\gamma\in \mathbb{R}$, the generalized extreme value distribution is denoted by $G_\gamma$. It is supported on $\Egamma$ and is be given by, for $x \in \Egamma$,
$$G_\gamma(x) =
\begin{cases}
\exp\left(-(1+\gamma x)^{-\frac1\gamma}\right), &\text{for $\gamma \neq 0$},\\
\exp\left(-\e^{-x}\right), &\text{for $\gamma = 0$}.
\end{cases}$$

\begin{definition} \label{def: domain of attr}
We say that a random variable $Y$ with distribution function $F$ is in the domain of attraction of an extreme value distribution $G_{\gamma}$ for some $\gamma\in \mathbb{R}$ (written as $D(G_\gamma)$) if there exists a positive valued function $a$ and a real valued function $b$ such that as $t\rightarrow \infty$, on $\Egamma$,
\begin{equation}\label{cevm:doa}
t\prob\left[Y>a(t)y+b(t)\right]=t \overline{F}(a(t)y+b(t))\rightarrow
\begin{cases}
      (1+\gamma y)^{-\frac{1}{\gamma}}, &\text{for $\gamma \neq 0$},\\
      \e^{-y}, &\text{for $\gamma = 0$}.
\end{cases}
\end{equation}
\end{definition}
When $\gamma\neq 0$, the domain of attraction condition is related to regular variation in the following way.

If $\gamma>0$, then $F\in D(G_{\gamma})$ if and only if $\overline{F}\in \text{RV}_{-1/\gamma}$.

If $\gamma<0$, then as $t\rightarrow\infty$, $b(t)\rightarrow b(\infty)<\infty$ and $F\in D(G_{\gamma})$ if and only if $\overline{F}(b(\infty)-1/\cdot)\in \text{RV}_{1/\gamma}.$ Note that in this case $b(\infty)$ becomes the upper end point of the distribution function~$F$.

\begin{definition} [Conditional extreme value model]\label{defn: cevm}
The real valued random vector $(X, Y)$ satisfies \textit{conditional extreme value model} (CEVM) if
\begin{enumerate}
\item \label{cevm:doa:def}The marginal distribution of $Y$ is in the domain of attraction of an extreme value distribution $G_{\gamma}$, that is, there exists a positive valued function $a$ and a real valued function $b$ such that~\eqref{cevm:doa} holds on $\Egamma$.

\item There exists a positive valued function $\alpha$ and a real valued function $\beta$ and a non-null Radon measure $\mu$ on Borel subsets of $(-\infty,\infty) \times \Egamma$ such that
    \begin{enumerate}[(2A)]
      \item \label{basic1} $t \prob \left[ \left(\frac{X-\beta(t)}{\alpha(t)}, \frac{Y-b(t)}{a(t)}\right) \in \cdot \right] \vaguec \mu(\cdot)$ on $[-\infty, \infty] \times \Ebar$, as $t\to\infty$, and
      \item \label{nondegenerate} for each $y \in \Egamma$, $\mu ((-\infty, x] \times (y, \infty))$ is a nondegenerate distribution function in $x$.
    \end{enumerate}

\item \label{cevm:probability} The function $H(x)=\mu((-\infty,x]\times(0,\infty))$ is a probability distribution.
\end{enumerate}
\end{definition}

If $(X,Y)$ satisfy Conditions~\eqref{cevm:doa:def}--\eqref{cevm:probability}, then we say $(X,Y)\in CEVM(\alpha,\beta; a, b;\mu)$ in  $[-\infty,\infty]\times\Ebar$.

Note that $H$ is a nondegenerate probability distribution function by Condition~\eqref{nondegenerate}. Also, Condition~\eqref{basic1} is equivalent to the convergence
$$t \prob [X \le \alpha(t) x + \beta(t), Y > a(t)y + b(t)] \to \mu ((-\infty, x] \times (y, \infty))$$ for all $y \in \Egamma$ and continuity points $(x, y)$ of the measure $\mu$. Note that, for $(X,Y)\in CEVM(\alpha,\beta;a,b;\mu)$ in  $[-\infty,\infty]\times\Ebar$ we have, as $t\rightarrow\infty$,
$$ \prob \left[ \left. \frac{X - \beta(t)}{\alpha(t)} \le x \right| Y > b(t) \right] \to H(x),$$ which motivates the name of the model.

Occasionally, we shall also be interested in pairs of random variables $(X,Y)$ which satisfy Conditions~\eqref{basic1} and~\eqref{nondegenerate}, without any reference to Conditions~\ref{cevm:doa:def} and~\ref{cevm:probability}. We shall then say that $(X,Y)$ satisfies Conditions~\eqref{basic1} and~\eqref{nondegenerate} with parameters $(\alpha,\beta; a, b;\mu)$ on  $E$, where $\alpha$ and $\beta$ will denote the scaling and centering of $X$, $a$ and $b$ will denote the scale and centering of $Y$ and $\mu$ will denote the nondegenerate limiting distribution and $E$ is the space on which the convergence takes place. In Definition~\ref{defn: cevm}, we have $E=[-\infty, \infty] \times \Ebar$.

It was shown in Proposition 1 of \cite{heffernan:resnick:2007} that there exists functions $\psi_1(\cdot)$, $\psi_2(\cdot)$ such that,
$$\lim_{t\rightarrow\infty}\frac{\alpha(tx)}{\alpha(t)}=\psi_1(x) \qquad \text{and} \qquad \lim_{t\rightarrow\infty}\frac{\beta(tx)-\beta(t)}{\alpha(t)}=\psi_2(x)$$
and the above convergence holds uniformly on compact subsets of $(0,\infty).$ Then, necessarily, we must have, for some $\rho \in \mathbb R$, $\psi_1(x)=x^{\rho}, x>0$. Also from Theorem B.2.1 of \cite{Haan:Ferreira:2006} it follows that either $\psi_2$ is $0$ or, for some $k\in\mathbb R$ and $x>0$,
$$\psi_2(x) =
\begin{cases}
  \frac k\rho (x^\rho - 1), &\text{when $\rho \neq 0$},\\
  k \log x, &\text{when $\rho = 0$}.
\end{cases}$$
Note that, we have $\alpha\in RV_\rho$.

\begin{definition} \label{def: mrv}
The pair of nonnegative random variables $(Z_1,Z_2)$ is said to be \textit{standard multivariate regularly varying} on  $[0,\infty]\times(0,\infty]$ if, as $t\rightarrow\infty$
$$ t\prob\left[\left(\frac{Z_1}{t},\frac{Z_2}{t}\right)\in \cdot\right]\vaguec\nu(\cdot) \text{ in } M_+([0,\infty]\times(0,\infty]).$$ 
\end{definition}  
In such cases we have $(Z_1,Z_2)\in CEVM(t,0;t,0;\nu)$ in $[0,\infty]\times(0,\infty]$.  The above convergence implies that $\nu(\cdot)$ is homogeneous of order $-1$, that is, $$\nu(c\Lambda)=c^{-1}\nu(\Lambda) \text{ for all  } c>0$$ where  $\Lambda$ is a Borel subset of $[0,\infty]\times(0,\infty]$. By homogeneity arguments it follows that for $r>0,$
\begin{align*}
&\nu\{(x,y)\in[0,\infty]\times(0,\infty]:x+y>r,\frac{x}{x+y}\in\Lambda\}\\
&\qquad =r^{-1}\nu\{(x,y)\in[0,\infty]\times(0,\infty]:x+y>1,\frac{x}{x+y}\in\Lambda\}\\
&\qquad=:r^{-1}S(\Lambda),
\end{align*}
where $S$ is a measure on $\{(x, y): x+y=1, 0\le x<1\}$. The measure $S$ is called the \textit{spectral measure} corresponding to $\nu(\cdot)$, while the measure $\nu$ is called the standardized measure. It was shown in \cite{heffernan:resnick:2007} that  whenever $(X,Y)\in CEVM(\alpha,\beta;a,b;\mu)$ in $[-\infty,\infty]\times\Ebar$ with $(\psi_1(x),\psi_2(x))\neq(1,0)$, we have the standardization  $(f_1(X),f_2(Y))\in CEVM(t,0;t,0;\nu)$ on the cone $[0,\infty]\times(0,\infty]$, for some monotone transformations $f_1$ and $f_2$.  \cite{Das:Resnick:2008} showed that this standardized measure $\nu$ is not a product measure. Throughout this article  we assume that $(\psi_1(x),\psi_2(x))\neq(1,0)$ and consider the product of $X$ and $Y$. We remark that although the model can be standardized in this case, the standardization does not help one to conclude about the behavior of $XY$.
\end{section}
\begin{section}{A brief overview of the results} \label{sec: overview}
In this section we give a brief description of the results. First note that  if $(X,Y)\in CEVM(\alpha,\beta;a,b;\mu)$ on  $[-\infty,\infty]\times\Ebar$, then  $\alpha\in RV_{\rho}$ and $a\in RV_\gamma$, where $\alpha$ and $a$ were the scalings for $X$ and $Y$ respectively. While $Y\in D(G_\gamma)$ necessarily holds, it  need not a priori follow that $X\in D(G_{\rho})$.  We classify the problem according to the parameters $\gamma$ and $\rho$. We break the problem into four cases depending on whether the parameters $\gamma$ and $\rho$ are positive or negative. In Section~\ref{Transformations} we show that, depending on the properties of the scaling and centering parameters, we can first reduce the basic convergence in conditional model to an equivalent convergence with the limiting measure satisfying nondegeneracy condition in an appropriate cone. The reduction of the basic convergence helps us to compute the convergence of the product with ease in Section~\ref{mainresults}.

{\bf Case I: $\rho$ and $\gamma$ positive:} This is an easier case and the behavior is quiet similar to the classical multivariate extreme theory.  In Theorem~\ref{case1[a]}, we show that under appropriate tail condition on $X$, the product $XY$ has regularly varying tail of index $-1/(\rho+\gamma)$. It is not assumed that $X\in D(G_{\rho})$, but in Section~\ref{section:remarks on moments} we show that the tail condition is satisfied when $X\in D(G_{\rho})$. It may happen that $X$ is in some other domain of attraction but still the tail condition holds. We also present a situation where the tail condition may fail.

In all the remaining cases, at least one of the indices $\rho$ and $\gamma$ will be negative. A negative value for $\gamma$ will require that the upper endpoint of $Y$ is indeed $b(\infty)$. However, as has been noted below, the same need not be true for $\rho$, $X$ and $\beta(\infty)$. Yet, we shall assume that whenever $\rho<0$, the upper endpoint of the support of $X$ is $\beta(\infty)$. Further, we shall assume at least one of the factors $X$ and $Y$ to be nonnegative. If both the factors take negative values, then, the product of the negative numbers being positive, their left tails will contribute to the right tail of the product as well. In fact, it will become important in that case to compare the relative heaviness of the two contributions. This can be easily done by breaking each random variable into its negative and positive parts. For the product of two negative parts, the relevant model should be built on $(-X, -Y)$ and the analysis becomes same after that. While these details increase the bookkeeping, they do not provide any new insight into the problem. So we shall refrain from considering the situations where both $X$ and $Y$ take negative values, except in Subcases II(b) and II(d) below, where both $X$ and $Y$ are nonpositive and we get some interesting result about the lower tail of the product $XY$ easily.

{\bf Case II: $\rho$ and $\gamma$ negative:} In Section~\ref{Transformations}, we first reduce the basic convergence to an equivalent convergence where regular variation can play an important role. In this case both $b(t)$ and $\beta(t)$ have finite limits $b(\infty)$ and $\beta(\infty)$ respectively. Since $Y\in D(G_{\gamma})$, $b(\infty)$ is the right end point of $Y$. However, $\beta(\infty)$ need not be the right end point of $X$ in general, yet throughout we shall assume it to be so. In Section~\ref{Transformations}, we reduce the conditional model $(X,Y)\in CEVM(\alpha,\beta;a,b;\mu)$ to $(\widetilde{X},\widetilde{Y})$ which satisfies Conditions~\eqref{basic1} and~\eqref{nondegenerate} with parameters $(\widetilde{\alpha},0;\widetilde{a},0;\nu)$ on $[0,\infty]\times(0,\infty]$, where,
\begin{equation}\label{eq:x-y-transformation}
 \widetilde{X}=\frac1{\beta(\infty)-X}\quad \text{and}\quad \widetilde{Y}=\frac1{b(\infty)-Y},
\end{equation}
and  $\widetilde{\alpha}$ and $\widetilde{a}$ are some appropriate scalings and $\nu$ is a transformed measure.
Regular variation at the right end point plays a crucial role during the determination the product behavior in this case. Depending on the right end point, we break the problem into few subcases which are interesting.

{\bf Subcase II(a): $\beta(\infty)$ and $b(\infty)$ positive:} If the right end points are positive, then, without loss of generality, we assume them to be $1$. In Theorem~\ref{main theorem negative}, we show that if $X$ and $Y$ both have positive right end points, then $(1-XY)^{-1}$ has regularly varying tail of index $-1/|\rho|$, under some further sufficient moment conditions. In Section~\ref{section:example}, we give an example where the moment condition fails, yet the product shows the tail behavior predicted by Theorem~\ref{main theorem negative}.

{\bf Subcase II(b):  $\beta(\infty)$ and $b(\infty)$  zero:} In Theorem~\ref{theorem:both zero} we show that if both right end points are zero, then the product convergence is a simple consequence of the result in Case I. In this case $(XY)^{-1}$ has regularly varying tail of index ${-1/(|\rho|+|\gamma|)}$.

{\bf Subcase II(c): $\beta(\infty)$ zero and $b(\infty)$ positive:} We show in Theorem~\ref{theorem: X zero} that if $Y$ is a nonnegative random variable having positive right end point, then under some appropriate moment conditions $-(XY)^{-1}$ has regularly varying tail of index ${-1/|\rho|}$.

{\bf Subcase II(d): $\beta(\infty)$ and $b(\infty)$  negative:} When both the right end points are negative, then, without loss of generality, we assume them to be $-1$. In Theorem \ref{theo: both right end point negative}, we show that $(XY-1)^{-1}$ has regularly varying tail of index ${-1/|\rho|}$.

There are a few more cases beyond the four subcases considered above, when both $\rho$ and $\gamma$ are negative. For example, consider the case when $Y$ has right end point zero and $X$ has positive right end point $\beta(\infty)$. By our discussion above, $X$ should have the support $[0,\beta(\infty)]$. Again, the product has right end point $0$ and the behavior of $X$ around zero becomes important. Thus, to get something interesting in this case one must have a conditional model which gives adequate information about the behavior of $X$ around the left end point. So it becomes natural to model $(-X, Y)$, which has already been considered in Subcase II(b). A similar situation occurs when $\beta(\infty)<0$ and $b(\infty)>0$. Here, again, the problem reduces to that in Subcase II(d) by modeling $(X, -Y)$.  We refer to Remark~\ref{rem: not done} for a discussion on this subcase.

{\bf Case III: $\rho$ positive and $\gamma$ negative:} In this case we assume $b(\infty)>0$ and also $\alpha(t)\sim 1/a(t)$ which implies that $\rho=-\gamma$. We  show in  Theorem~\ref{theorem case 3} that $XY$ has regularly varying tail of index ${-1/|\gamma|}$.

{\bf Case IV: $\rho$ negative and  $\gamma$ positive:} In Theorem~\ref{case 4} we show that $XY$ has regularly varying tail of index ${-1/\gamma}$.

Finally we end this section by summarizing the results in a tabular form:

\begin{center}
\begin{table}[h]
\caption{Behavior  of  products}	\label{tab:BehaviorOfProducts}
	\begin{tabular}{|c|c|c| c| c| c|}
	\hline \hline
	Index of $\alpha$     & Index of $a$     &Theorem number  & Nature     & Regular variation    \\[0.5ex]
	\hline
	$\rho>0$ & $\gamma >0$ & Theorem~\ref{case1[a]} & $XY$ &$\text{RV}_{-{1}/{(\gamma+\rho)}}$ \\
	\hline
	$\stackrel{\displaystyle \rho>0}{\alpha\sim\frac{1}{a}}$ & $\gamma=-\rho<0$  &  Theorem~\ref{theorem case 3} & $XY$ & $\text{RV}_{-{1}/{|\gamma|}}$ \\
	\hline
	$\rho<0$ & $\gamma<\rho$ &  Theorem~\ref{main theorem negative}& $(\beta(\infty)b(\infty)-XY)^{-1}$ & $\text{RV}_{-{1}/{|\rho|}}$ \\
	\hline
	$\stackrel{\displaystyle\rho<0}{\alpha\sim a}$ & $\gamma=\rho$ & Theorem~\ref{main theorem negative} &  $(\beta(\infty)b(\infty)-XY)^{-1}$ &$\text{RV}_{-{1}/{|\rho|}}$ \\
	\hline	
$\rho<0$& $\gamma<0$  & Theorem~\ref{theorem:both zero} & $(XY)^{-1}$ & $\text{RV}_{-{1}/(|\gamma|+|\rho|)}$ \\
\hline
$\rho<0$& $\gamma<0$  & Theorem~\ref{theorem: X zero}& $-(XY)^{-1}$ & $\text{RV}_{-{1}/{|\rho|}}$ \\
\hline
$\rho<0$& $\gamma>0$ & Theorem~\ref{case 4} & XY& $\text{RV}_{1/\gamma}$\\
\hline
\hline
\end{tabular}	
\end{table}
\end{center}

\end{section}

\begin{section}{Some transformations of CEVM according to parameters $\gamma$ and $\rho$} \label{Transformations}
In this section we reduce the basic convergence in Condition~\eqref{basic1} to an equivalent convergence in some appropriate subspace of $\mathbb R^2$ to facilitate  our calculations of  product of two variables following conditional extreme value model. We now discuss the four cases considered in Section~\ref{sec: overview}.

{\bf Case I: $\rho$ and $\gamma$ positive:} In this case we assume $Y$ is nonnegative. Let $(X,Y)\in CEVM(\alpha,\beta;a,b;\mu)$ on $[-\infty,\infty]\times \Ebar$. Now by the domain of attraction condition~\eqref{cevm:doa:def} and  Corollary 1.2.4 of \cite{Haan:Ferreira:2006}, we have,  $b(t)\sim a(t)/\gamma$, as $t\to\infty$. Also, from Theorem 3.1.12 (a),(c) of \cite{bingham:goldie:teugels:1987} it follows that
  $$\lim_{t\rightarrow\infty}\frac{\beta(t)}{\alpha(t)}=
  \begin{cases}
  0 &\text{when $\psi_2=0$} \\
  \frac{1}{\rho} &\text{when $\psi_2\neq 0$}.
  \end{cases}$$
Now using the above conditions and translating $X$ and $Y$ coordinates we get that $(X,Y)$ satisfies Conditions~\eqref{basic1} and~\eqref{nondegenerate} with parameters $(\alpha,0;a,0;\nu)$ on $D:=[-\infty,\infty]\times(0,\infty]$, for some nondegenerate measure $\nu$ which is obtained from $\mu$ by translations on both axes. So in Theorem~\ref{case1[a]} which deals with product $XY$ in Case I, we assume that $(X,Y)$ satisfies Conditions~\eqref{basic1} and~\eqref{nondegenerate} with parameters $(\alpha,0;a,0;\nu)$ on $D=[-\infty,\infty]\times(0,\infty]$ for some nondegenerate Radon measure $\nu$.

{\bf Case II: $\rho$ and $\gamma$ negative:}  Recall that in this case $\Ebar=\left(-\infty,\frac{1}{|\gamma|}\right]$.  Since $Y\in D(G_{\gamma})$ with $\gamma<0$ , it follows from Lemma 1.2.9 of \cite{Haan:Ferreira:2006} that  $\lim_{t\rightarrow\infty}b(t)=:b(\infty)$ exists and is finite and as $t\rightarrow\infty$ we have, $$\frac{b(\infty)-b(t)}{a(t)}\rightarrow\frac{1}{|\gamma|}.$$ Moreover $b(\infty)$ turns out to be the right end point of $Y$. Hence in this case, without loss of generality we take $a(t)=|\gamma|(b(\infty)-b(t))$ and it easily follows that, for $y>0$,
$$\lim_{t\rightarrow\infty}t\prob\left[\frac{\widetilde Y}{a(t)^{-1}}>y\right]=  y^{\frac{1}{\gamma}},$$ where $\widetilde Y$ is defined in~\eqref{eq:x-y-transformation}.  Now observe that $(X,Y)\in CEVM(\alpha,\beta;a,b;\mu)$ on $[-\infty,\infty]\times \Ebar$ gives $(X,\widetilde Y)$ which satisfies Conditions~\eqref{basic1} and~\eqref{nondegenerate} with parameters $(\alpha,\beta;1/a(t),0;\mu_2)$ on $D$, where, $$\mu_2 ([-\infty,x]\times(y,\infty])=\mu( [-\infty,x]\times(\frac{1}{|\gamma|}-\frac{1}{y},\infty]).$$

Now since $\rho<0$, we get by  Theorem B.22 of \cite{Haan:Ferreira:2006} that $\lim_{t\rightarrow\infty}\beta(t)=\beta(\infty)$ exists and is finite. It may happen that $X$ has a different  right end point than $\beta(\infty)$,  but we assume $\beta(\infty)$ to be its right end point to avoid complications. In  $X$ coordinate we can do a similar transformation as the $Y$ variable, to get
\begin{align*}
K_t(\alpha(t)x+\beta(t))&:=\prob\left[  \left. \frac{X-\beta(t)}{\alpha(t)}\leq x \right|\frac{(b(\infty)-Y)^{-1}}{a(t)^{-1}}>y\right]\\
&\rightarrow y^{-\frac1\gamma}\mu_2([-\infty,x]\times(y,\infty])=:K(x).
\end{align*}
When $\psi_2\neq0$, we have, as $t\to\infty$,
\begin{equation}\label{eq:beta:negative:conv}
\frac{\beta(\infty)-\beta(t)}{\alpha(t)}\rightarrow\frac{1}{|\rho|}.\end{equation}

Now by convergence of types theorem, as $t\to\infty$, we have, $$K_t(|\rho|(\beta(\infty)-\beta(t))x+\beta(t))\rightarrow K(x).$$ Define,
\begin{equation}\label{eq:negative:alpha:a}
\widetilde{\alpha}(t)=
\begin{cases}
\frac1{|\rho|(\beta(\infty)-\beta(t))} &\text{when~} \psi_2\neq0\\
\frac1{\alpha(t)} &\text{when~} \psi_2=0\\
\end{cases}
\qquad \text{and} \qquad
\widetilde{a}(t)=\frac1{a(t)}.
\end{equation}
Using~\eqref{eq:negative:alpha:a} we get, as $t\rightarrow\infty$,
\begin{equation}
\label{basic2a}
t\prob\left[\frac{\widetilde X}{\widetilde{\alpha}(t)}\leq x,\frac{\widetilde Y}{\widetilde{a}(t)}>y\right]\to
\begin{cases}
\mu_2([-\infty,-\frac{1}{x}+\frac{1}{|\rho|}]\times(y,\infty]) &\text{for }\psi_2\neq 0\\
 \mu_2([-\infty,-\frac{1}{x}]\times(y,\infty]) &\text{for }\psi_2=0.\\
\end{cases}
\end{equation}
In Section~\ref{mainresults}, we deal with Case II by breaking it up into different subcases as pointed out in Section~\ref{sec: overview}. So, in Theorems~\ref{main theorem negative}--\ref{theo: both right end point negative}, we assume that $(\widetilde X,\widetilde Y)$  satisfy Conditions~\eqref{basic1} and~\eqref{nondegenerate} with parameters $(\widetilde{\alpha},0;\widetilde{a},0;\nu)$ on $[0,\infty]\times(0,\infty]$, for some nondegenerate Radon measure $\nu$.

{\bf Case III:  $\rho$ positive and $\gamma$ negative: }
Since $Y\in D(G_{\gamma})$, we can do a transformation similar to that in Case II. So in this case  $(X,\widetilde Y)$  satisfies Conditions~\eqref{basic1} and~\eqref{nondegenerate} with parameters $(\alpha,\beta;\widetilde{a},0;\mu_3)$ on $[-\infty,\infty]\times(0,\infty]$, for some nondegenerate measure $\mu_3$.

Now since $\rho>0$, we can do a translation in the first coordinate to get $(X,\widetilde Y)$  which satisfies Conditions~\eqref{basic1} and~\eqref{nondegenerate} with parameters $(\alpha,0;\widetilde{a},0;\nu)$ on $[-\infty,\infty]\times(0,\infty]$ for some nondegenerate measure $\nu$. In Theorem~\ref{theorem case 3} we deal with product behavior in this case.

{\bf Case IV:  $\rho$ negative and $\gamma$ positive: }
We assume that $\beta(\infty)$, the right end point of X is positive. Now, as in Case II, we use~\eqref{eq:beta:negative:conv} to get the following convergence for $x\geq 0$ and $y>0,$
\begin{align*}
   t\prob\left[\frac{(\beta(\infty)-X)}{\alpha(t)}\leq x,\frac{Y}{a(t)}>y\right]&=t\prob\left[\frac{X-\beta(t)}{\alpha(t)}\geq -x+\frac{\beta(\infty)-\beta(t)}{\alpha(t)},\frac{Y}{a(t)}>y\right]\\
&\rightarrow\mu\left([-x+\frac{1}{|\gamma|},\infty]\times (y,\infty]\right) \text{ as } t\rightarrow\infty.
\end{align*}
So in Theorem~\ref{case 4}, which derives the product behavior in Case IV, we assume  $(\beta(\infty)-X, Y)$  satisfies Conditions~\eqref{basic1} and~\eqref{nondegenerate} with parameters $(\alpha,0;a,0;\nu)$ on $[0,\infty]\times(0,\infty]$, for some nondegenerate Radon measure $\nu$.

We thus observe that if $(X, Y)\in CEVM(\alpha, \beta; a, b; \mu)$, then in Cases I, II, III and IV respectively, $(X, Y)$, $(\widetilde X, \widetilde Y)$, $(X, \widetilde Y)$ and $(\beta(\infty)-X, Y)$ satisfy Conditions~\eqref{basic1} and~\eqref{nondegenerate} with some positive scaling parameters, zero centering parameters and a nondegenerate limiting Radon measure on $D=[-\infty, \infty]\times (0, \infty]$. In future sections, whenever we refer to Conditions~\eqref{basic1} and~\eqref{nondegenerate} with respect to the transformed variables alone without any reference to the CEVM model for the original pair $(X, Y)$, we shall denote, by an abuse of notation, the limiting Radon measure for the transformed random variables as $\mu$ as well.
\end{section}

\begin{section}{Behavior of the product under conditional model}
\label{mainresults}
Now we study the product behavior when $(X,Y)\in CEVM(\alpha,\beta;a,b;\mu)$. In all the cases we assume Conditions~\eqref{basic1} and~\eqref{nondegenerate} on the suitably transformed versions of $(X, Y)$, so that centering is not required.

{\bf Case I: $\rho$ and $\gamma$ positive:} We begin with the case where both $\rho$ and $\gamma$ are positive. 
\begin{theorem}\label{case1[a]} Let $\rho>0,\gamma>0$ and $Y$ be a nonnegative random variable. Assume $(X,Y)$ satisfies Conditions~\eqref{basic1} and~\eqref{nondegenerate} with parameters $(\alpha,0;a,0;\mu)$ on $D:=[-\infty,\infty]\times(0,\infty]$. Also assume
  \begin{equation}
  \label{moment condition}
   \lim_{\epsilon\downarrow0}\limsup_{t\rightarrow\infty}t\prob\left[\frac{|X|}{\alpha(t)}>\frac{z}{\epsilon}\right]=0.
  \end{equation}
   Then, $XY$ has regularly varying tail of index $-{1/(\gamma+\rho)}$ and $t\prob\left[XY/{\left(\alpha(t)a(t)\right)}\in\cdot\right]$ converges vaguely to some nondegenerate Radon measure on $[-\infty,\infty]\setminus \{0\}$.
\end{theorem}
\begin{proof}
For $\epsilon>0$ and $z>0$ observe that the set, $$A_{\epsilon,z}=\{(x,y)\in D :xy>z, y>\epsilon \}$$ is a relatively compact set in $D$ and $\mu$ is a Radon measure. Note that,
\begin{multline}
t\prob\left[\left(\frac{XY}{\alpha(t)a(t)},\frac{Y}{a(t)}\right)\in A_{\epsilon,z}\right]\le t\prob\left[\frac{XY}{\alpha(t)a(t)}>z\right] \\ \le  t\prob\left[\left(\frac{XY}{\alpha(t)a(t)},\frac{Y}{a(t)}\right)\in A_{\epsilon,z}\right] + t\prob\left[\frac{XY}{\alpha(t)a(t)}>z,\frac{Y}{a(t)}\leq \epsilon\right]
\end{multline}
First letting $t\to\infty$ and then $\epsilon\downarrow0$ through a sequence such that $A_{\epsilon,z}$ is a $\mu$-continuity set, the left side and the first term on the right side converge to $\mu\{(x,y)\in D: xy>z\}$. The second term on the right side is negligible by the assumed tail condition~\eqref{moment condition}. Combining all, we have the required result when $z>0$.
By similar arguments one can show the above convergence on the sets of the form $(-\infty,-z)$ with $z>0$, also.
\end{proof}

\noindent{\bf Spectral form for product:}
For simplicity let us assume that $X$ is nonnegative as well. Then the vague convergence in $M_+(D)$ can be thought of as vague convergence in $M_+([0,\infty]\times(0,\infty]).$  By Theorem \ref{case1[a]} we have,
$$\lim_{t\rightarrow\infty}t\prob\left[\frac{XY}{\alpha(t)a(t)}>z\right]= \mu\{(x,y)\in [0,\infty]\times(0,\infty] : xy > z\}.$$
Since $\rho>0$ and $\gamma>0$, it is known that there exists $\overline{\alpha}(t)\sim \alpha(t)$ and $\overline{a}(t)\sim a(t)$ such that they are eventually differentiable and strictly increasing. Also $\frac{\overline{\alpha}(tx)}{\alpha(t)}\rightarrow x^{\rho}$ and $\frac{\overline{a}(tx)}{a(t)}\rightarrow x^{\gamma}$ as $t\rightarrow\infty.$  Recall that the (left continuous) inverse of nondecreasing function $f$ is defined as $$f^{\leftarrow}(y)=\inf \{s: \ f(s)\ge y\ \}.$$
Hence,
\begin{align*}
   t\prob\left[\frac{\overline{\alpha}^{\leftarrow}(X)}{t}\leq x, \frac{\overline{a}^{\leftarrow}(Y)}{t}>y\right]&=t\prob\left[\frac{X}{\alpha(t)}\leq \frac{\overline{\alpha}(tx)}{\alpha(t)},\frac{Y}{a(t)}>\frac{\overline{a}(ty)}{a(t)}\right]\\
   &\rightarrow \mu([0,x^{\rho}]\times(y^{\gamma},\infty])\\
   &=\mu T_1^{-1}([0,x]\times(y,\infty]),
 \end{align*}
where $T_1(x,y)=(x^{1/\rho},y^{1/\gamma})$.
Let $S$ be the spectral measure for the standardized pair $(\overline{\alpha}^{\leftarrow}(X),\overline{a}^{\leftarrow}(Y))$ corresponding to $\mu T_1^{-1}$. Then,
  \begin{align*}
  &\mu\{(x,y)\in[0,\infty]\times(0,\infty]: xy > z\}\\
  =&\mu T_1^{-1}\{(x,y)\in [0,\infty]\times(0,\infty] : x^{\rho}y^{\gamma} > z\}\\
  =&\int_{\omega\in[0,1)}\int_{r^{\rho+\gamma}\omega^{\rho}(1-\omega)^{\gamma}>z} r^{-2}dr S(d\omega)\\
  =&\int_{\omega\in[0,1)}\int_{r>\frac{z^{\frac{1}{\rho+\gamma}}}{(\omega^{\rho}(1-\omega)^{\gamma})^{\frac{1}{\rho+\gamma}}}}r^{-2}drS(d\omega)\\
  =&z^{-\frac{1}{\rho+\gamma}}\int_{\omega\in[0,1)}\omega ^{\frac{\rho}{\rho+\gamma}}(1-\omega)^{\frac{\gamma}{\rho+\gamma}}S(d\omega).\\
  \end{align*}
So finally we have,
\begin{align*}\lim_{t\rightarrow\infty}t\prob\left[\frac{XY}{\alpha(t)a(t)}>z\right]&= \mu\{(x,y)\in [0,\infty]\times(0,\infty] : xy > z\}\\
&=z^{-\frac{1}{\rho+\gamma}}\int_{\omega\in[0,1)}\omega ^{\frac{\rho}{\rho+\gamma}}(1-\omega)^{\frac{\gamma}{\rho+\gamma}}S(d\omega).\\
\end{align*}

{\bf Case II: $\rho$ and $\gamma$ negative:} As has been already pointed out, $\beta(\infty)$ need not be the right endpoint $X$. However, we shall assume it to be so. The tail behavior of $XY$ strongly depends on the right end points of $X$ and $Y$. There are several possibilities which may arise, but it may not always be possible to predict the tail behavior of $XY$ in all the cases. We shall deal with few interesting cases. See Section~\ref{sec: overview} for a discussion in this regard. Regarding one of the cases left out, see Remark~\ref{rem: not done}. Recall that, from~\eqref{eq:x-y-transformation} we have,
\begin{equation}\label{transformation}
\frac{1}{\beta(\infty)b(\infty)-XY}=\frac{\widetilde{X}\widetilde{Y}}{\beta(\infty)\widetilde{X}+b(\infty)\widetilde{Y}-1}
\end{equation}

{\bf Subcase II(a): $\beta(\infty)$ and $b(\infty)$ positive:} After scaling $X$ and $Y$ suitably, without loss of generality, we can assume that $\beta(\infty)=1=b(\infty).$
\begin{theorem}\label{main theorem negative}
Suppose $X$ and $Y$ are nonnegative and $(\widetilde X,\widetilde Y)$ satisfies Conditions~\eqref{basic1} and~\eqref{nondegenerate} with parameters $(\widetilde{\alpha},0;\widetilde{a},0;\mu)$ on $[0,\infty]\times(0,\infty]$. Assume $\E\left[\widetilde{X}^{{1}/{|\rho|}+\delta}\right]<\infty$ for some $\delta>0$ and  either  $\gamma < \rho $ or  $\widetilde{\alpha}(t)/\widetilde{a}(t)$ remains bounded. Then $(1-XY)^{-1}$ has regularly varying tail of index ${-1/|\rho|}$ and, as $t\to\infty$, $t\prob\left[{(1-XY)^{-1}}/{\widetilde{\alpha}(t)}\in\cdot \right]$ converges vaguely to some nondegenerate Radon measure on $(0,\infty]$.
\end{theorem}

We start with a technical lemma.
\begin{lemma}\label{important lemma}
For $0\leq t_1\leq t_2\leq\infty$ and $z>0$, we denote the set
\begin{equation}\label{eq:set:V}
V_{[t_1,t_2],z}=\{(x,y)\in [0,\infty]\times(0,\infty]:x\in[t_1,t_2], y>z\}.
\end{equation} 
Suppose that $\{(Z_{1t},Z_{2t})\}$ is a sequence of pairs of nonnegative random variables and there exists a Radon measure $\nu(\cdot)$ in $[0,\infty]\times(0,\infty]$  such that they satisfy the following two conditions.\\
Condition A: Whenever $V_{[t_1,t_2],z}$ is a $\nu$-continuity set, we have, as $t\rightarrow\infty$,
$$t\prob\left[(Z_{1t},Z_{2t})\in V_{[t_1,t_2],z}\right]\rightarrow \nu(V_{[t_1,t_2],z}).$$
Condition B: For any  $z_0\in(0,\infty)$ we have as $t\rightarrow\infty,$
$$t\prob\left[(Z_{1t},Z_{2t})\in V_{[0,\infty],z_0}\right]\rightarrow f(z_0)\in(0,\infty).$$
Then, as $t\to\infty$,
$$t\prob\left[(Z_{1t},Z_{2t})\in\cdot\right]\vaguec\nu(\cdot)$$ in $M_+([0,\infty]\times(0,\infty])$
\end{lemma}
\begin{proof}
Fix a $z_0\in(0,\infty)$ and define the following probability measures on $[0,\infty]\times(z_0,\infty)$:
$$Q_t(\cdot)=\frac{t\prob\left[(Z_{1t},Z_{2t})\in \cdot\right]}{t\prob\left[(Z_{1t},Z_{2t})\in V_{[0,\infty],z_0}\right]} \text{ and } Q(\cdot)=\frac{\nu(\cdot)}{f(z_0)}.$$ From condition A it follows that $$Q_t(V_{[t_1,t_2],z})\rightarrow Q(V_{[t_1,t_2],z})$$  as $t\rightarrow\infty$, whenever $V_{[t_1,t_2],z}$ is $\nu$-continuity set. Now following the arguments in the proof of Theorem 2.1 of \cite{Maulik:Resnick:Rootzen:2002}, it follows that $Q_t$ converges weakly to $Q$ on $[0,\infty]\times(z_0,\infty)$. Since a Borel set with boundary having zero $Q$ measure is equivalent to having measure zero with respect to the measure $\nu$ we have that $t \prob\left[(Z_{1t},Z_{2t})\in B\right]\rightarrow\nu(B)$ for any Borel set $B$ having boundary with zero $\nu$ measure.

Let $K$ be a $\nu$-continuity set as well as a relatively compact set in $[0,\infty]\times(0,\infty]$. Then there exists $z_0>0$ such that $K\subset[0,\infty]\times(z_0,\infty]$. Then $K$ is Borel in $[0,\infty]\times(z_0,\infty)$ and also a $\nu$-continuity set. Hence we have,
$$t\prob\left[(Z_{1t},Z_{2t})\in K\right]=t\prob\left[(Z_{1t},Z_{2t})\in K\right]\rightarrow \nu(K).$$ This shows that $t\prob\left[(Z_{1t},Z_{2t})\in \cdot\right]$ vaguely converges to $\nu$ on $[0,\infty]\times(0,\infty]$.
\end{proof}

From~\eqref{transformation}, the behavior of $XY$ will be determined by the pair $(\widetilde X \widetilde Y, \widetilde X + \widetilde Y - 1)$. So we next prove a result about the joint convergence of product and the sum.
\begin{lemma}\label{lemma:joint: prod:sum}
Let $\gamma<0,\rho<0$ and $(\widetilde X,\widetilde Y)$  satisfies Conditions~\eqref{basic1} and~\eqref{nondegenerate} with parameters $(\widetilde{\alpha},0;\widetilde{a},0;\mu)$. If  $\E\left[\widetilde{X}^{{1}/{|\gamma|}+\delta}\right]<\infty$ for some $\delta>0,$ then~$(\widetilde X\widetilde Y,\widetilde X+\widetilde Y-1)$ also satisfies Conditions~\eqref{basic1} and~\eqref{nondegenerate} with parameters $(\widetilde{\alpha}\widetilde{a},0;\widetilde{a},0;\mu T_2^{-1})$ on  $[0,\infty]\times(0,\infty]$ where $T_2(x,y)=(xy,y)$.
\end{lemma}
\begin{proof}
First observe that from the compactification arguments used in Lemma~\ref{important lemma} and the basic convergence in Condition~\eqref{basic1} satisfied by the pair $(\widetilde X,\widetilde Y)$, it follows that
\begin{equation} \label{T_2 convergence}
t\prob\left[\left(\frac{\widetilde{X}\widetilde{Y}}{\widetilde{\alpha}(t)\widetilde{a}(t)},\frac{\widetilde{Y}}{\widetilde{a}(t)}\right)\in\cdot\right]\vaguec \mu T_2^{-1}(\cdot)\quad \text{in}\quad  M_+([0,\infty]\times(0,\infty]).
\end{equation}
Let $0\leq t_1\leq t_2\leq \infty$ and $z>0$. Assume that $\mu T_2^{-1}(\partial V_{[t_1,t_2],z})=0$, where $V_{[t_1,t_2],z}$ is defined in~\eqref{eq:set:V}. Since $X$ is nonnegative, $\widetilde{X}$ is greater than or equal to $1$ and hence for a lower bound we get,
\begin{multline}
\liminf_{t\rightarrow\infty}t\prob\left[\left(\frac{\widetilde{X}\widetilde{Y}}{\widetilde{\alpha}(t)\widetilde{a}(t)}, \frac{\widetilde{X}+\widetilde{Y}-1}{\widetilde{a}(t)}\right)\in V_{[t_1,t_2],z}\right]\\ \geq\lim_{t\rightarrow\infty}t\prob\left[\frac{\widetilde{X}\widetilde{Y}}{\widetilde{\alpha}(t)\widetilde{a}(t)}\in[t_1,t_2],\frac{\widetilde{Y}}{\widetilde{a}(t)}>z\right]
=\mu T_2^{-1}(V_{[t_1,t_2],z}).\label{eq:negative:lower bound}
\end{multline}
For the upper bound, choose $0<\epsilon<z$, such that $1/|\widetilde a(t)| < \epsilon/2$ (since $\widetilde{a}(t)\in RV_{-\gamma}$) and $\mu T_2^{-1}\left(\partial V_{[t_1,t_2],z-\epsilon}\right)=0$. Then,
\begin{align*}
&t\prob\left[\left(\frac{\widetilde{X}\widetilde{Y}}{\widetilde{\alpha}(t)\widetilde{a}(t)},\frac{\widetilde{X}+\widetilde{Y}-1}{\widetilde{a}(t)}\right)\in V_{[t_1,t_2],z}\right]\\
\leq &t\prob\left[\left(\frac{\widetilde{X}\widetilde{Y}}{\widetilde{\alpha}(t)\widetilde{a}(t)},\frac{\widetilde{X}+\widetilde{Y}-1}{\widetilde{a}(t)}\right)\in V_{[t_1,t_2],z},\frac{\widetilde{X}}{\widetilde{a}(t)}\leq\frac{\epsilon}2\right]\\
& +t\prob\left[\left(\frac{\widetilde{X}\widetilde{Y}}{\widetilde{\alpha}(t)\widetilde{a}(t)},\frac{\widetilde{X}+\widetilde{Y}-1}{\widetilde{a}(t)}\right)\in V_{[t_1,t_2],z},\frac{\widetilde{X}}{\widetilde{a}(t)}>\frac{\epsilon}2\right]\\
\leq &t\prob \left[\frac{\widetilde{X}\widetilde{Y}}{\widetilde{\alpha}(t)\widetilde{a}(t)}\in[t_1,t_2],\frac{\widetilde{Y}}{\widetilde{a}(t)}>z-\epsilon\right]
+t\prob\left[\frac{\widetilde{X}}{\widetilde{a}(t)}>\frac{\epsilon}2\right]\\
\leq &t\prob\left[\left(\frac{\widetilde{X}\widetilde{Y}}{\widetilde{\alpha}(t)\widetilde{a}(t)},\frac{\widetilde{Y}}{\widetilde{a}(t)}\right)\in V_{[t_1,t_2],z-\epsilon}\right] + 2^{1/|\gamma|+\delta} t \frac{\E\left[\widetilde{X}^{{1}/{|\gamma|}+\delta}\right]}{(\widetilde{a}(t)\epsilon)^{\frac{1}{|\gamma|}+\delta}}.
\end{align*}
The first term converges to $\mu T_2^{-1}(V_{[t_1,t_2],z-\epsilon})$, while the second sum converges to zero, since $\widetilde{a}(t)\in RV_{-\gamma}$. Now letting $\epsilon\to0$ satisfying the defining conditions, we obtain the upper bound, which is same as the lower bound~\eqref{eq:negative:lower bound}. Thus we get that,
$$\lim_{t\rightarrow\infty}t\prob\left[(\frac{\widetilde{X}\widetilde{Y}}{\widetilde{\alpha}(t)\widetilde{a}(t)},\frac{\widetilde{X}+\widetilde{Y}}{\widetilde{a}(t)})\in V_{[t_1,t_2],z}\right]= \mu T_2^{-1}(V_{[t_1,t_2],z}).$$
Hence Condition A of Lemma~\ref{important lemma} is satisfied.

Now if we fix $z_0\in(0,\infty)$  and let $\epsilon>0$ satisfy the conditions as in the upper bound above, then
\begin{align*} t\prob\left[\left(\frac{\widetilde{X}\widetilde{Y}}{\widetilde{\alpha}(t)\widetilde{a}(t)},\frac{\widetilde{X}+\widetilde{Y}-1}{\widetilde{a}(t)}\right)\in V_{[0,\infty],z_0}\right]&= t\prob\left[\frac{\widetilde{X}+\widetilde{Y}-1}{\widetilde{a}(t)}>z_0\right]\\
&\leq t\prob\left[\frac{\widetilde{Y}}{\widetilde{a}(t)}>z_0-\epsilon \right]+t\prob\left[\frac{\widetilde{X}}{\widetilde{a}(t)}>\frac{\epsilon}2 \right]\\
&\rightarrow (z_0 -\epsilon)^{\frac{1}{\gamma}}.
\end{align*}
Hence the upper bound for the required limit in Condition B of Lemma~\ref{important lemma} follows by letting $\epsilon\rightarrow 0.$ The lower bound easily follows from the domain of attraction condition on $Y$ and the fact that $\widetilde X\ge 1$. So Condition B is also satisfied by the pair~$(\frac{\widetilde{X}\widetilde{Y}}{\widetilde{\alpha}(t)\widetilde{a}(t)},\frac{\widetilde{X}+\widetilde{Y}-1}{\widetilde{a}(t)})$ and hence the result follows from Lemma \ref{important lemma}.
\end{proof}

\begin{proof}[Proof of Theorem~\ref{main theorem negative}]
Denote $W'=\widetilde{X}\widetilde{Y}$, $W''=\widetilde{X}+\widetilde{Y}-1$ and note that $(1-XY)^{-1}=W'/W''$. So from previous lemma it follows  that,
$$t\prob\left[\left(\frac{W'}{\widetilde{\alpha}(t)\widetilde{a}(t)},\frac{W''}{\widetilde{a}(t)}\right)\in \cdot\right]\vaguec \mu T_2^{-1}(\cdot)\qquad \mbox{as~} t\rightarrow\infty.$$
Let $w\in(0,\infty)$ and $\epsilon>0$ and consider the set,
\begin{equation}\label{eq: neg: A set:descrip}
B_{w,\epsilon}=\{(x,y)\in [0,\infty]\times(0,\infty]: x>yw, y>\epsilon\}.\end{equation}
Then, for $\epsilon>0$, we have,
$$t\prob\left[\frac{W'}{W''}\frac{1}{\widetilde{\alpha}(t)}> w \right] = t \prob\left[\left(\frac{W'}{\widetilde{\alpha}(t)\widetilde{a}(t)}, \frac{W''}{\widetilde{a}(t)}\right)\in B_{w,\epsilon}\right] + t \prob\left[\frac{W'}{W''}\frac{1}{\widetilde{\alpha}(t)}> w, \frac{W''}{\widetilde{a}(t)}\leq \epsilon \right].$$
Since the set is bounded away from both the axes, the first sum converges to $\mu T_2^{-1}(B_{w,\epsilon})$ by the vague convergence of $(W', W'')$. Now since $\widetilde{X}\geq0$ we have $\widetilde{X}+\widetilde{Y}-1\geq \widetilde{Y}-1$, and hence for large $t$ we get,
\begin{align*}
&t\prob\left[\frac{(1-XY)^{-1}}{\widetilde{\alpha}(t)}>w, \frac{\widetilde{Y}-1}{\widetilde{a}(t)}\leq \epsilon_k\right]\\
\leq &t\prob\left[ XY>1-\frac{1}{w\widetilde{\alpha}(t)}, Y\leq 1-\frac{1}{\widetilde{a}(t)\epsilon_k +1} \right]\\
&\leq t\prob\left[\widetilde{X}>\frac{1-\frac{1}{\widetilde{a}(t)\epsilon_k+1}}{\frac{1}{w\widetilde{\alpha}(t)}-\frac{1}{(\widetilde{a}(t)\epsilon_k+1)}}\right]\\
&\leq C(t,w,k)\E\left[\widetilde{X}^{\frac1{|\rho|}+\delta}\right],
\end{align*}
where,
$$C(t,w,k)= \frac{t}{(\widetilde{\alpha}(t))^{\frac{1}{|\rho|}+\delta}}\left({1-\frac{1}{\widetilde{a}(t)\epsilon_k+1}}\right)^{-(\frac{1}{|\rho|}+\delta)} \left(\frac{1}{w}-\frac{1}{(\frac{\widetilde{a}(t)}{\widetilde{\alpha}(t)}\epsilon_k+\frac{1}{\widetilde{\alpha}(t)})}\right)^{\frac{1}{|\rho|}+\delta},$$
which goes to zero as $t\to\infty$, since $\widetilde{\alpha}(t)\in RV_{-\rho}$ and $\widetilde{\alpha}(t)/\widetilde{a}(t)$ remains bounded.
\end{proof}
\begin{remark}
Theorem~\ref{main theorem negative} requires that $(1/|\rho|+\delta)$-th moment of $\widetilde X$ is finite. However, this condition is not necessary. In the final section we give an example where this moment condition is not satisfied but we still obtain the tail behavior of the product.
\end{remark}

{\bf Subcase II(b): $\beta(\infty)=0$ and $ b(\infty)=0$:} In this case, both $X$ and $Y$ are nonpositive, but the product $XY$ is nonnegative. Thus, the right tail behavior of $XY$ will be controlled by the left tail behaviors of $X$ and $Y$, which we cannot control much using CEVM. However, CEVM gives some information about the left tail behavior of $XY$ at $0$, which we summarize below. 

Note that in this case, from~\eqref{eq:x-y-transformation} and~\eqref{eq:negative:alpha:a}, we have $\widetilde{X}=-1/{X}$, $\widetilde{Y}=-{1}/{Y}$ and $\widetilde{\alpha}(t)=-{1}/{(|\rho|\beta(t))}$, $\widetilde{a}(t)={1}/{a(t)}$. From Theorem~\ref{case1[a]}, the behavior of the product $XY$ around zero, or equivalently the behavior of the reciprocal of the product $\widetilde X\widetilde Y$ around infinity, follows immediately.
\begin{theorem}\label{theorem:both zero}
If $\rho<0$, $\gamma<0$ and $(\widetilde X,\widetilde Y)$  satisfies Conditions~\eqref{basic1} and~\eqref{nondegenerate} with parameters $(\widetilde{\alpha},0;\widetilde{a},0;\mu)$.  Also suppose, $$\lim_{\epsilon\downarrow0}\limsup_{t\rightarrow\infty}t\prob\left[\frac{\widetilde{X}}{\widetilde{\alpha}(t)}>\frac{z}{\epsilon}\right]=0$$  Then  $(XY)^{-1}$ has regularly varying tail with index $-1/(|\gamma|+|\rho|)$ and as $t\rightarrow\infty$, $t\prob\left[{(XY)^{-1}}/({\widetilde{\alpha}(t)\widetilde{a}(t)})\in\cdot\right]$ converge to some nondegenerate Radon measure on $(0,\infty]$.
\end{theorem}

{\bf Subcase II(c): $\beta(\infty)=0$ and $b(\infty)=1$:} Now note that from~\eqref{transformation} we have, $$-\frac{1}{XY}=\frac{\widetilde{X}\widetilde{Y}}{\widetilde{Y}-1}.$$
\begin{theorem}\label{theorem: X zero}
Suppose $Y$ is nonnegative and $(\widetilde X,\widetilde Y)$  satisfies Conditions~\eqref{basic1} and~\eqref{nondegenerate} with parameters $(\widetilde{\alpha},0;\widetilde{a},0;\mu)$. If $\E\left[\widetilde{X}^{{1}/{|\rho|}+\delta}\right]<\infty$ for some $\delta>0$, then  $-(XY)^{-1}$ has regularly varying tail of index $-{1/|\rho|}$ and as $t\rightarrow\infty$, $t\prob\left[{-(XY)^{-1}}/{\widetilde{\alpha}(t)}\in\cdot\right]$ converge vaguely to some nondegenerate measure on $(0,\infty]$.
\end{theorem}
\begin{proof}
Under the given hypothesis, and the fact that $\widetilde{a}(t)\to \infty$, it can be shown by arguments similar to Lemma~\ref{lemma:joint: prod:sum} that,
$ t\prob\left[\left({\widetilde{X}\widetilde{Y}}/({\widetilde{\alpha}(t)\widetilde{a}(t)}), ({\widetilde{Y}-1})/{\widetilde{\alpha}(t)}\right)\in\cdot\right]$ converges vaguely to some nondegenerate Radon measure in $M_+([0,\infty]\times(0,\infty]$.

Next for $z>0$ and $\epsilon>0$, so that the set $B_{z,\epsilon}$ as in~\eqref{eq: neg: A set:descrip} is a continuity set of the limit measure. Note that
\begin{multline*}
 t\prob\left[\left(\frac{\widetilde{X}\widetilde{Y}}{\widetilde{\alpha}(t)\widetilde{a}(t)}, \frac{\widetilde{Y}-1}{\widetilde{\alpha}(t)}\right)\in B_{z,\epsilon}\right]\leq t\prob\left[\frac{-(XY)^{-1}}{\widetilde{\alpha}(t)}>z\right]\\ \leq t\prob\left[\left(\frac{\widetilde{X}\widetilde{Y}}{\widetilde{\alpha}(t)\widetilde{a}(t)}, \frac{\widetilde{Y}-1}{\widetilde{\alpha}(t)}\right)\in B_{z,\epsilon}\right]+t\prob\left[\frac{\widetilde{X}\widetilde{Y}}{\widetilde{\alpha}(t)(\widetilde{Y}-1)} >z,\frac{\widetilde{Y}-1}{\widetilde{a}(t)}\leq\epsilon \right].
\end{multline*}

The term on the left side and the first term on the right side converge as $t\to\infty$ due to the vague convergence mentioned at the beginning of the proof. Letting $\epsilon\downarrow0$ appropriately, we get the required limit. For the second term on the right side observe that,
 \begin{align*}
t\prob\left[\frac{\widetilde{X}\widetilde{Y}}{\widetilde{\alpha}(t)(\widetilde{Y}-1)} >z,\frac{\widetilde{Y}-1}{\widetilde{a}(t)}\leq\epsilon \right]&
\leq t\prob\left[\frac{-(XY)^{-1}}{\widetilde{\alpha}(t)}>z,\widetilde{Y}\leq \widetilde{a}(t)\epsilon+1\right]\\
&\leq t\prob\left[\frac{-(XY)^{-1}}{\widetilde{\alpha}(t)}>z,1-Y\geq\frac{1}{\widetilde{a}(t)\epsilon+1}\right]\\
&\leq t\prob\left[\frac{XY}{\alpha(t)}>-\frac{1}{z},Y\leq 1-\frac{1}{\widetilde{a}(t)\epsilon+1}\right]\\
&\leq t\prob\left[-X< \frac{\alpha(t)/z}{1-\frac{1}{\widetilde{a}(t)\epsilon+1}}\right]=t\prob\left[\widetilde{X}> \frac{1-\frac{1}{\widetilde{a}(t)\epsilon+1}}{\alpha(t)/z}\right]\\
&\leq t\E\left[\widetilde{X}^{{1}/{|\rho|}+\delta}\right]\left(\frac{\alpha(t)/z}{1-\frac{1}{\widetilde{a}(t)\epsilon+1}}\right)^{\frac{1}{|\rho|}+\delta}.
\end{align*}
The last expression tends to zero as $\widetilde{a}(t)\rightarrow\infty$ and $\alpha\in RV_\rho$.
\end{proof}

{\bf Subcase II(d): $\beta(\infty)$ and $b(\infty)$ negative:} As in Subcase II(a), after suitable scaling, without loss of generality, we can assume that $\beta(\infty)=b(\infty)=-1$.
Again, from~\eqref{transformation}, we have
$$\frac1{XY-1}=\frac{\widetilde{X}\widetilde{Y}}{\widetilde{X}+\widetilde{Y}+1}$$
and to get the behavior of the product around $1$ we first need to derive the joint convergence of $(\widetilde{X}\widetilde{Y}, \widetilde{X}+\widetilde{Y}+1)$. Using an argument very similar to Theorem~\ref{main theorem negative}, we immediately obtain the following result.
\begin{theorem} \label{theo: both right end point negative}
Let $(\widetilde X,\widetilde Y)$  satisfy Conditions~\eqref{basic1} and~\eqref{nondegenerate} with parameters $(\widetilde{\alpha},0;\widetilde{a},0;\mu)$ and $\E\left[\widetilde{X}^{\frac{1}{|\rho|}+\delta}\right]<\infty$ for some $\delta>0$. If either $\gamma < \rho $ or $\widetilde{\alpha}(t)/\widetilde{a}(t)$ remains bounded, then $(XY-1)^{-1}$ has regularly varying tail of index ${-1/|\rho|}$ and as $t\rightarrow\infty$,
$t\prob\left[{(XY-1)^{-1}}/{\widetilde{\alpha}(t)}\in\cdot \right]$ converges vaguely to some nondegenerate Radon measure on $(0,\infty]$.
\end{theorem}

\begin{remark} \label{rem: not done}
The other case when $\beta(\infty)=1$ and $b(\infty)=0$ is not easy to derive from the information about the conditional extreme value model. In this case the right endpoint of the product is zero and behavior of $X$ around zero seems to be important. But the conditional model gives us the regular variation behavior around one and not around zero.
\end{remark}

{\bf Case III: $\rho$ positive and $\gamma$ negative:} In this case we shall assume that $X$ is nonnegative and the upper endpoint of $Y$, $b(\infty)$ is positive. If $b(\infty)\leq0$, then the behavior of $X$ around its lower endpoint will play a crucial role in the behavior of the product $XY$, which becomes negative. However, the behavior of $X$ around its lower endpoint is not controlled by the conditional model and we are unable to conclude about the product behavior when $b(\infty)\leq 0$. So we only consider the case $b(\infty)>0$. We also make the assumption that $\alpha(t)\sim\widetilde a(t)$, which requires that $\rho=|\gamma|$.
\begin{theorem}\label{theorem case 3}
Let $X$ be nonnegative and $Y$ have upper endpoint $b(\infty)>0$. Assume that $(X,\widetilde Y)$  satisfies Conditions~\eqref{basic1} and~\eqref{nondegenerate} with parameters$(\alpha,0;\widetilde{a},0;\mu)$ with $\alpha(t)\sim {1}/{a(t)}=\widetilde{a}(t)$ and $\E\left[X^{{1}/{|\gamma|}+\delta}\right]<\infty$ for some $\delta>0$, then $XY$ has regularly varying tail of index ${-1/|\gamma|}$ and as $t\rightarrow\infty$, we have $t\prob\left[{XY}/\widetilde a(t) \in\cdot\right]$ converges vaguely to a nondegenerate Radon measure on $(0,\infty]$.
\end{theorem}
\begin{proof} 
As $\alpha(t)\sim{\widetilde a(t)}$, convergence of types allows us to change $\alpha(t)$ to $\widetilde{a}(t)$ and hence $(X,\widetilde Y)$ satisfy Conditions~\eqref{basic1} and~\eqref{nondegenerate} with parameters $(\widetilde{a},0;\widetilde{a},0;\mu)$.  Using the fact that $\widetilde{Y}=(b(\infty)-Y)^{-1}$, we have  $$XY=X\frac{b(\infty)\widetilde{Y}-1}{\widetilde{Y}}.$$ Using arguments similar to Proposition~4 of \cite{heffernan:resnick:2007}, it can be shown that as $t\rightarrow\infty$,
$$t\prob\left[\left(\frac{X}{\widetilde{Y}},\frac{b(\infty)\widetilde{Y}}{\widetilde{a}(t)}\right)\in\cdot\right]\vaguec \mu T_3^{-1}(\cdot),$$ in $M_+([0,\infty]\times(0,\infty])$, where $T_3(x,y)=(x/y,b(\infty)y)$. Since $\widetilde{a}(t)\rightarrow\infty$, we further have, as $t\rightarrow\infty,$
$$t\prob\left[\left(\frac{X}{\widetilde{Y}},\frac{b(\infty)\widetilde{Y}-1}{\widetilde{a}(t)}\right)\in\cdot\right]\vaguec \mu T_3^{-1}(\cdot).$$
Now applying the map $T_2(x,y)=(xy,y)$ to the above vague convergence and using compactification arguments similar to that in the proof of Theorem~2.1 of \cite{Maulik:Resnick:Rootzen:2002}, we have,
$$t\prob\left[\left(\frac{X}{\widetilde{Y}} \frac{b(\infty)\widetilde{Y}-1}{\widetilde{a}(t)}, \frac{b(\infty)\widetilde{Y}-1}{\widetilde{a}(t)}\right)\in\cdot\right]\vaguec \mu T_3^{-1}T_2^{-1}(\cdot).$$
Recalling the facts that $XY = X (b(\infty) \widetilde Y - 1)/\widetilde Y$ and $\widetilde a(t)\to\infty$ and reversing the arguments in the second coordinate, we have,
\begin{equation} \label{vague 3}
t\prob\left[\left(\frac{XY}{\widetilde{a}(t)},\frac{\widetilde{Y}}{\widetilde{a}(t)}\right)\in\cdot\right]\vaguec \mu T_3^{-1}T_2^{-1}{\widetilde{T}_3}^{-1}(\cdot),
\end{equation}
where $\widetilde{T}_3(x,y)=(x,\frac{y}{b(\infty)}).$   If $\epsilon>0$ and $z>0$ we have following series of inequalities,
\begin{align*}
   t\prob\left[\frac{XY}{\widetilde{a}(t)}>z,\frac{\widetilde{Y}}{\widetilde{a}(t)}\leq\epsilon\right]&
   =t\prob\left[\frac{XY}{\widetilde{a}(t)}>z,\frac{1}{b(\infty)-Y}\leq \widetilde{a}(t)\epsilon\right]\\
   &=t\prob\left[\frac{XY}{\widetilde{a}(t)}>z,Y\leq b(\infty)-\frac{1}{\widetilde{a}(t)\epsilon}\right]\\
   &\leq t\prob\left[\frac{X}{\widetilde{a}(t)}>z\left(b(\infty)-\frac{1}{\widetilde{a}(t)\epsilon}\right)^{-1}\right]\\
   &\leq t\frac{\E\left[X^{\frac{1}{|\gamma|}+\delta}\right]}{\left(\widetilde{a}(t)z\right)^{\frac{1}{|\gamma|}+\delta}}
   \left(b(\infty)-\frac{1}{\widetilde{a}(t)\epsilon}\right)^{\frac{1}{|\gamma|}+\delta} \to 0,
\end{align*}
since $\widetilde a \in RV_{|\gamma|}$.

Now observe that
\begin{equation*}
t \prob\left[\frac{XY}{\widetilde a(t)}>z, \frac{\widetilde Y}{\widetilde a(t)}>\epsilon\right] \le t \prob\left[\frac{XY}{\widetilde a(t)}>z\right]
\le  t \prob\left[\frac{XY}{\widetilde a(t)}>z, \frac{\widetilde Y}{\widetilde a(t)}>\epsilon\right] + t \prob\left[\frac{XY}{\widetilde a(t)}>z, \frac{\widetilde Y}{\widetilde a(t)}\le\epsilon\right].\end{equation*}
The last term on the right side is negligible by the previous argument and the left side and the first term on the right side converge due to the vague convergence in~\eqref{vague 3}. The result then follows by letting $\epsilon\to 0$.
\end{proof}

{\bf Case IV: $\rho$ negative and $\gamma$ positive:}
In this case we shall assume that $Y$ is nonnegative, $\beta(\infty)>0$ and $\beta(\infty)$ is the upper endpoint of $X$. Arguing as in Case III, we neglect the possibility that $\beta(\infty)\leq0$. Thus we further assume $0\leq X\leq \beta(\infty)$. Since $X$ becomes bounded, the product of $XY$ inherits its behavior from the tail behavior of $Y$.
\begin{theorem}\label{case 4}
Assume that both $X$ and $Y$ are nonnegative random variables with $\beta(\infty)$ being the upper endpoint of $X$. Let $(\beta(\infty)-X, Y)$   satisfies Conditions~\eqref{basic1} and~\eqref{nondegenerate} with parameters $(\alpha,0;a,0;\mu)$ on $[0,\infty]\times(0,\infty]$, for some nondegenerate Radon measure $\mu$. Then $XY$ has regularly varying tail of index ${-1/\gamma}$ and for all $z>0$, we have
$$t\prob\left[\frac{XY}{a(t)}>z\right]\rightarrow z^{-\frac1\gamma}\beta(\infty)^{\frac1\gamma}, \text{ as } t\rightarrow \infty.$$
\end{theorem}
\begin{proof}
First we prove the upper bound which, in fact, does not use the conditional model. Fix $\epsilon>0$ and observe that
\begin{align*}
 t\prob\left[\frac{XY}{a(t)}>z\right]&=t\prob\left[\frac{XY}{a(t)}>z,\epsilon<X<\beta(\infty)\right]+t\prob\left[\frac{XY}{a(t)}>z,X\leq \epsilon\right]\\
&\leq t\prob\left[\frac{Y}{a(t)}>\frac{z}{\beta(\infty)}\right] +t\prob\left[\frac{Y}{a(t)}>\frac{z}{\epsilon}\right]\\
&\rightarrow \left(\frac{z}{\beta(\infty)}\right)^{-\frac{1}{\gamma}}+\left(\frac{z}{\epsilon}\right)^{-\frac{1}{\gamma}} \text{ as } t\rightarrow\infty.
\end{align*}
The second term goes to zero as $\epsilon\downarrow 0$. 

To prove the lower bound we use the basic convergence in Condition~\eqref{basic1} for the pair $(\beta(\infty)-X, Y)$. Before we show the lower bound, first observe that by arguments similar to the proof of Theorem 2.1 of \cite{Maulik:Resnick:Rootzen:2002} we have,
\begin{equation}\label{basic4b}
 t\prob\left[\left(\frac{(\beta(\infty)-X)Y}{\alpha(t)a(t)},\frac{Y}{a(t)}\right)\in\cdot\right]\vaguec \mu T_2^{-1}(\cdot) \text{ as } t\rightarrow\infty
\end{equation}
in $M_+([0,\infty]\times(0,\infty])$, where recall that $T_2(x,y)=(xy,y).$ Now to show the lower bound, first fix a large $M>0$ and $\epsilon>0$ such that $\alpha(t)<\epsilon$ for large $t$ (recall that $\alpha\in RV_{-\rho}$ and $\alpha(t)\rightarrow0$ in this case). Note that
\begin{align*}
 t\prob\left[\frac{XY}{a(t)}>z\right]&=t\prob\left[\frac{(X-\beta(\infty))Y}{\alpha(t)a(t)}\alpha(t)+\frac{\beta(\infty)Y}{a(t)}>z\right]\\
&\geq t\prob\left[\frac{(\beta(\infty)-X)Y}{\alpha(t)a(t)}\leq M, \frac{\beta(\infty)Y}{a(t)}>z+M\alpha(t)\right]\\
&\geq t\prob\left[\frac{(\beta(\infty)-X)Y}{\alpha(t)a(t)}\leq M, \frac{\beta(\infty)Y}{a(t)}>z+M\epsilon\right]\\
&\to \mu T_2^{-1}\left([0,M]\times \left(\frac{z+M\epsilon}{\beta(\infty)},\infty\right]\right),
\end{align*}
using~\eqref{basic4b}. First letting $\epsilon\rightarrow0$ and then letting $M\rightarrow\infty$, so that $$\mu T_2^{-1}\left(\partial \left([0, M]\times (z/\beta(\infty), \infty]\right)\right)=0,$$ we get that
$$\liminf_{t\rightarrow\infty}t\prob[\frac{XY}{a(t)}>z]\geq \mu T_2^{-1}([0,\infty)\times (z/\beta(\infty),\infty])=z^{-\frac1\gamma}\beta(\infty)^{\frac1\gamma}.$$
\end{proof}
\end{section}

\begin{section}{Some remarks on the tail condition~\eqref{moment condition}}\label{section:remarks on moments}
We have already noted that the compactification arguments in Section~\ref{mainresults} require some conditions on the tails of associated random variables. Except in Theorem~\ref{case1[a]}, they have been replaced by some moment conditions, using Markov inequality. The tail condition~\eqref{moment condition} in Theorem~\ref{case1[a]} can also be replaced by the following moment condition:

If for some $\delta>0$, we have $\E[|X|^{{1}/{\rho}+\delta}]<\infty$, then~\eqref{moment condition} holds, as $\alpha\in RV_{\rho}$.

In general, if $(X,Y)$ follows CEVM model, it need not be true that $X\in D(G_\rho)$. However, if $X\in D(G_{\rho})$  with scaling and centering functions $\alpha(t)$ and $\beta(t)$ and $X\geq0$, then the moment condition \eqref{moment condition}. In fact,
$$\limsup_{t\rightarrow\infty}t\prob\left[\frac{X}{\alpha(t)}>\frac{z}{\epsilon}\right]$$ is a constant multiple of $({z}/{\epsilon})^{-{1}/{\rho}}$, which goes to zero as $\epsilon\rightarrow 0.$

The tail condition~\eqref{moment condition} continues to hold in certain other cases as well. Suppose that $X\geq0$ and $X\in D(G_{\lambda})$ with scaling and centering $A(t)$ and $B(t)$ and $\lambda<\rho$. In this case, $\alpha\in RV_\rho$ and $A\in RV_\lambda$. Thus, $\alpha(t)/A(t)\to\infty$. Hence, for any $\epsilon>0$, we have
$$t\prob\left[\frac{X}{\alpha(t)}>\frac{z}{\epsilon}\right]=t\prob\left[\frac{X}{A(t)}>\frac{z\alpha(t)}{\epsilon A(t)}\right]$$
is of order of $({\alpha(t)}/{A(t)})^{-{1}/{\lambda}}$, which goes to zero and~\eqref{moment condition} holds. However, it would be interesting to see the effect of $A$ and $B$ as scaling and centering in CEVM model. Since, $\alpha$ is of an order higher than $A$, the limit, as expected, becomes degenerate.

\begin{proposition}
Let the pair $(X,Y)\in CEVM(\alpha,\beta;a,b;\mu)$  with $\rho>0$ and $\gamma>0.$ Assume  $X\in D(G_{\lambda})$ with  $\rho>\lambda$ and centering and scaling as $A(t)$ and $B(t)$. If 
$$\lim_{t\rightarrow\infty}t\prob\left[\frac{X-B(t)}{A(t)}\leq x,\frac{Y-b(t)}{a(t)}>y\right]$$ exists for all continuity points $(x,y)\in \mathbb{R}\times \mathbb{E}^{(\gamma)}$, then, for any fixed $y\in \mathbb{E}^{(\gamma)}$, as $x$ varies in $\mathbb R$, the limit measure assigns same values to the sets of the form $[-\infty, x]\times (y,\infty]$.
\end{proposition}
\begin{proof} 
Observe that
$$\frac{A(t)}{\alpha(t)}x+\frac{B(t)-\beta(t)}{\alpha(t)}=\frac{A(t)}{\alpha(t)}\left(x+\frac{B(t)}{A(t)}\right)-\frac{\beta(t)}{\alpha(t)}\rightarrow 0\times(x+\frac{1}{\lambda})-\frac{1}{\rho}.$$
Therefore we have,
\begin{multline*}
\lim_{t\rightarrow\infty} t\prob\left[\frac{X-B(t)}{A(t)}\leq x,\frac{Y-b(t)}{a(t)}>y\right]\\
=\lim_{t\rightarrow\infty} t\prob\left[\frac{X-\beta(t)}{\alpha(t)}\leq \frac{A(t)}{\alpha(t)}x+\frac{B(t)-\beta(t)}{\alpha(t)},\frac{Y-b(t)}{a(t)}>y\right]
=\mu( [-\infty,-\frac{1}{\rho}]\times(y,\infty]),
\end{multline*}
which is independent of $x$.
\end{proof}

We would also like to consider the remaining case, namely, when $X\ge0$ and $X\in D(G_\lambda)$ with scaling $A$ and $\lambda>\rho$. Clearly, we have $\alpha(t)/A(t)\to 0$. Thus, for any $\epsilon>0$, we have 
$$t\prob\left[\frac{X}{\alpha(t)}>\frac{z}{\epsilon}\right]=t\prob\left[\frac{X}{A(t)}>\frac{z\alpha(t)}{\epsilon A(t)}\right]\to\infty$$ Hence~\eqref{moment condition} cannot hold and Theorem~\ref{case1[a]} is of no use. The next result show that in this case we have multivariate extreme value model with the limiting measure being concentrated on the axes, which gives asymptotic independence.
\begin{theorem}\label{theo:hidden}
Let $(X,Y)\in CEVM(\alpha,0;a,0;\mu)$  on the cone $[0,\infty]\times(0,\infty]$ and $X\in D(G_\lambda)$ with $\lambda>\rho$ and scaling $A(t)\in RV_{-\lambda}$.  Also assume that $X\in D(G_\lambda)$ with $\lambda>\rho$ and ${A(t)}/{\alpha(t)}\rightarrow \infty$. Then
\begin{equation}\label{hidden2}
t\prob\left[\left(\frac{X}{A(t)},\frac{Y}{a(t)} \right)\in\cdot\right]
\end{equation}
converges vaguely to a nondegenerate Radon measure on $[0,\infty]^2 \setminus\{(0,0)\}$, which is concentrated on the axes.
\end{theorem}
\begin{proof}
We first show the vague convergence on $[0,\infty]^2 \setminus\{(0,0)\}$. Let $x>0,y>0$. Then,
\begin{multline*}t\prob\left[ \left(\frac{X}{A(t)},\frac{Y}{a(t)}\right)\in ([0,x]\times [0,y])^c\right]\\ = t\prob\left[\frac{X}{A(t)}>x\right]+t\prob\left[\frac{Y}{a(t)}>y\right]-t\prob\left[\frac{X}{A(t)}>x,\frac{Y}{a(t)}>y\right].\end{multline*}
The first two terms converge due to the domain of attraction conditions on $X$ and $Y$. For the last term, by the CEVM conditions and the fact that $A(t)/\alpha(t)\to\infty$, we have
\begin{equation} \label{eq:hidden}
t\prob\left[\frac{X}{A(t)}>x,\frac{Y}{a(t)}>y\right] = t\prob\left[\frac{X}{\alpha(t)}>x\frac{A(t)}{\alpha(t)},\frac{Y}{a(t)}>y\right]\to 0.
\end{equation}
This establishes~\eqref{hidden2}. However, using $x=y$ in~\eqref{eq:hidden}, we find that the limit measure does not put any mass on $(0,\infty]^2$.
\end{proof}
\begin{remark}
Since, we have asymptotic independence, several different behaviors for the product are possible, as it has been illustrated in \cite{Maulik:Resnick:Rootzen:2002}.

While we have established asymptotic independence on the larger cone $[0,\infty]^2\setminus\{0\}$ in Theorem~\ref{theo:hidden}, CEVM gives another nondegenerate limit $\mu$ on the smaller cone $[0,\infty]\times(0,\infty]$. Thus, $(X,Y)$ exhibits hidden regular variation, as described in \cite{Resnick:2002, Maulik:Resnick:2004}.
\end{remark}
\end{section}

\begin{section}{Example}\label{section:example}
We now consider the moment condition in Theorem~\ref{main theorem negative}. We show that the condition is not necessary by providing an example, where the condition fails, but we explicitly calculate the tail behavior of $XY$. 

Let $X$ and $Z$ be two independent random variables, where $X$ follows Beta distribution with parameters $1$ and $a$ and $Z$ is supported on $[0,1]$ and is in $D(G_{-1/b})$, for some $a>0$ and $b>0$. Thus, we have $\prob[X>x]=(1-x)^a$ and $\prob[Z>1-\frac{1}{x}]=x^{-b}L(x)$ for some slowly varying function $L$. Let $G$ denote the distribution function of the random variable $Y=X\wedge Z$. Then $\overline{G}(x)=(1-x)^{a+b}L(\frac{1}{1-x})$ and hence $Y\in D(G_{-1/(a+b)}).$ Clearly, for $\widetilde X= 1/(1-X)$, we have $\E[\widetilde X^{a+b}]=\infty$ and the moment condition in Theorem~\ref{main theorem negative} fails for $\rho=-1/(a+b)$.

We further define  
$$\widetilde{a}(t)= \frac{1}{1-\left({1}/{\overline{G}}\right)^{\leftarrow}(t)}.$$ Since, $Y\in D(G_{-1/(a+b)})$, we have from Corollary 1.2.4 of \cite{Haan:Ferreira:2006}, \begin{equation}\label{example relation}
 \frac{{\widetilde{a}(t)}^{a+b}}{L(\widetilde{a}(t))}\sim t.
\end{equation}
Then, for $x>y>0$, we have,
\begin{align*}
t \prob \left[ \frac{\widetilde X}{\widetilde a(t)} \le x, \frac{\widetilde Y}{\widetilde a(t)} > y\right]
= &t\prob \left[\widetilde{a}(t)(X-1)\leq -\frac1x, \widetilde{a}(t)(Y-1)> -\frac1y\right]\\
= &t\prob \left[1-\frac{1}{y\widetilde{a}(t)}<X\leq 1-\frac{1}{x\widetilde{a}(t)}, Z>1-\frac{1}{y\widetilde{a}(t)}\right]\\
= &\frac{t}{\widetilde{a}(t)^{a+b}}[y^{-a}-x^{-a}]y^{-b}L(y{\widetilde{a}(t)})\\
\sim &\frac{t L(\widetilde{a}(t))}{\widetilde{a}(t)^{a+b}}[y^{-a}-x^{-a}]y^{-b}
\sim [y^{-a}-x^{-a}]y^{-b},
\end{align*}
where we use~\eqref{example relation} in the last step. Thus, $(\widetilde X, \widetilde Y)$ satisfies Conditions~\eqref{basic1} and~\eqref{nondegenerate} with $\rho=\gamma=-1/(a+b)$.

Finally, we directly calculate the asymptotic tail behavior of the product. For simplicity of the notations, for $y>0$, let us denote $c(t)=1-{1}/(\widetilde{a}(t)y)$. Since $\widetilde{a}\in RV_{1/(a+b)}$, as $t\to\infty$, we have $c(t)\to 1$. Also, $\widetilde a(t) = 1/((1-c(t))y)$. Then,
\begin{align*}
&t\prob\left[XY>1-\frac{1}{\widetilde{a}(t)y}\right]
= a\cdot t \int\limits_{\sqrt{c(t)}}^{1} \prob \left[Z>\frac{c(t)}{s} \right](1-s)^{a-1}ds\\
\sim &a\cdot {\widetilde{a}}^{\leftarrow} \left(\frac{1}{(1-c(t))y}\right) \int\limits_{\sqrt{c(t)}}^{1}\prob \left[Z>\frac{c}{s}\right] (1-s)^{a-1}ds\\
\sim &a\cdot {\widetilde{a}}^{\leftarrow} \left(\frac{1}{(1-c(t))y}\right) \int\limits_{\sqrt{c(t)}}^{1} \left(1-\frac{c}{s}\right)^{b} (1-s)^{a-1} L\left(\frac{1}{1-c(t)/s}\right)ds
\intertext{substituting $s=1-z(1-c(t))$,}
= &a (1-c(t))^{a+b} {\widetilde{a}}^{\leftarrow} \left(\frac{1}{(1-c(t))y}\right) \int\limits_{0}^{\frac{1}{1+\sqrt{c(t)}}} \frac{(1-z)^{b}z^{a-1}}{(1-(1-c(t))z)^{b}}  L\left(1+\frac{c(t)}{(1-c(t))(1-z)}\right)dz\\
\sim &a (1-c(t))^{a+b} {\widetilde{a}}^{\leftarrow} \left(\frac{1}{(1-c(t))y}\right) L\left( \frac{c(t)}{1-c(t)}\right) \int\limits_{0}^{1/2}(1-z)^b z^{a-1}dz,
\end{align*}
where, in the last step, we use Dominated Convergence Theorem and the facts that $c(t)\to 1$ and 
$$L\left(1+\frac{c(t)}{(1-c(t))(1-z)}\right) \sim L\left(\frac{c(t)}{(1-c(t))(1-z)}\right) \sim L\left( \frac{c(t)}{1-c(t)}\right)$$
uniformly on bounded intervals of $z$.
Finally, using~\eqref{example relation}, definition of $c(t)$,  to get,
$$t\prob\left[\frac{(1-XY)^{-1}}{\widetilde{a}(t)}>y\right]\rightarrow a \cdot y^{-(a+b)}\int\limits_{0}^{1/2}(1-z)^b z^{a-1}dz.$$
\end{section}

\end{document}